\input amstex
\documentstyle{amsppt}

\magnification=\magstep1

\pagewidth{5.5truein}
\hcorrection{0.5truein}
\pageheight{8.5truein}
\vcorrection{0.25truein}

\NoPageNumbers
\NoBlackBoxes
\TagsAsMath
\TagsOnLeft

\def\NN{{\Bbb N}}

\def\Aus{{\bf 1}}
\def\Auslarge{{\bf 2}}
\def\ChaI{{\bf 3}}
\def\ChaII{{\bf 4}}
\def\Cohkap{{\bf 5}}
\def\FulRei{{\bf 6}}
\def\Gri{{\bf 7}}
\def\GruJen{{\bf 8}}
\def\Her{{\bf 9}}
\def\Herz{{\bf 10}}
\def\Hui{{\bf 11}}
\def\Huibiel{{\bf 12}}
\def\Hull{{\bf 13}}
\def\Koe{{\bf 14}}
\def\Lenok{{\bf 15}}
\def\RinTac{{\bf 16}}
\def\Simcounter{{\bf 17}}
\def\Simpad{{\bf 18}}
\def\War{{\bf 19}}
\def\Zim{{\bf 20}}
\def\Zimakad{{\bf 21}}

\topmatter

\title  Direct products of modules and the pure
semisimplicity conjecture 
\endtitle

\rightheadtext{products and pure semisimplicity}

\author Birge Huisgen-Zimmermann and Frank Okoh \endauthor

\thanks The research of the first author was partially
supported by a grant from the National Science Foundation.
\endthanks

\address Department of Mathematics, University of
California, Santa Barbara, CA 93106\endaddress

\email birge\@math.ucsb.edu\endemail

\address Department of Mathematics, Wayne State University,
Detroit, MI 48202 
\endaddress

\email okoh\@math.wayne.edu \endemail

\dedicatory  Professor Helmut Lenzing anlaesslich seines
60.~Geburtstags gewidmet
\enddedicatory

\endtopmatter

\document

\magnification=\magstep1

Whenever global decomposition problems involving {\it all}
modules over a ring are at stake, the `resistance' of
infinite direct products to being decomposed into
manageable direct summands plays a crucial role.  This
phenomenon already surfaced in work of Koethe and
Cohen-Kaplansky \cite{\Koe, \Cohkap} and in Chase's 
landmark paper \cite{\ChaI}, as well as in numerous
articles that have appeared since, e.g., in work of 
Warfield, Griffith, Auslander, Gruson-Jensen, W.
Zimmermann, and the first author (see, e.g., \cite{\War,
\Gri, \Auslarge, \GruJen,
\Zim, \Hui}).  While, in the meantime, major headway has
been made towards identifying and analyzing the rings whose
module categories enjoy the most useful decomposition
properties  --  such as decomposability of all left modules
into indecomposable direct summands or, equivalently, into
finitely generated components  --  it is still not clear
whether these demands on the base ring are left-right
symmetric.  The problem is all the more tantalizing as the
two-sided requirement is known to coincide with finite
representation type (\cite{\Aus},
\cite{\RinTac}, \cite{\FulRei});  resolving it thus amounts
to locating the missing piece towards a link between finite
and infinite dimensional representation theory.  We refer to
\cite{\Simpad} and
\cite{\Huibiel} for more detail on the current status of
the problem.  The question on which we focus here was
raised in its present form by the second author, in
pursuance of problems he had jointly tackled with Lenzing
\cite{\Lenok}.  It is also representative of the hurdles
that remain in the way of understanding the symmetry issue
we sketched.  Indeed, a positive answer for all two-sided
artinian rings would settle symmetry in the positive.  

\definition{Question} For which rings  R  does the
following hold:  If
$(M_n)_{n \in
\NN}$ is any family of pairwise non-isomorphic finitely
generated, indecomposable left
$R$-modules, then the direct product $\prod_{n \in \NN} M_n$
fails to be a direct sum of finitely generated submodules? 
\enddefinition 

We point out that a
negative answer to the symmetry question for global
decompositions is expected
\cite{\Simcounter}, whence, a fortiori, existence of artinian
counterexamples to our product condition appears likely.

The purpose of this note is to exhibit several important
instances in which the product condition is satisfied; we
hope that our sketch of the base line will trigger further
work in this direction.  In particular, we will see that
all Artin algebras are among the rings satisfying the
condition, as well as all commutative noetherian domains of
Krull dimension 1.  With a modest amount of new input, one
can glean these results from the existing literature. 
While, in complete absence of finiteness conditions for
$R$, the outcome is negative  --  just let $R$ be an
infinite direct product of local rings $R_n$ and choose
$M_n = R_n$  --  we conjecture that our product condition
is satisfied by all noetherian PI rings.  In the noetherian
case, decomposability of
$\prod_{n \in
\NN} M_n$ into finitely generated submodules of course
implies decomposability into indecomposable summands.  In view
of our example, we will
significantly increase the class of rings falling on the
positive side by weakening the above product condition, so as to
only demand the following: There is no infinite family of
pairwise non-isomorphic finitely generated and indecomposable
left $R$-modules whose direct product is a direct sum
of finitely generated indecomposable components.

We start by making the connection of our question with
global decomposition problems more precise.  It is
well-known that the rings all of whose left modules are
direct sums of finitely generated (or indecomposable)
submodules are precisely those over which all pure
inclusions split \cite{\Zimakad, \Hui}.  Their name
reflects this behavior: they have been dubbed {\it left
pure semisimple rings}.  In \cite{\ChaI}, Chase observed
that they are necessarily left artinian.  As already
mentioned, it is known that left and right pure
semisimplicity is tantamount to finite representation type,
but so far, no examples are available indicating that
one-sided pure semisimplicity does not automatically carry
over to the other side.  In the presence of mild
commutativity conditions, the property is in fact known to
be left-right symmetric; namely, for Artin algebras
\cite{\Auslarge}, and, more generally, for Artin PI rings
and rings permitting a (Morita) self-duality \cite{\Her,
\Herz}.  

\proclaim{Observation 1}  If all two-sided artinian rings
satisfy the above product condition, then each left pure
semisimple ring has finite representation type.
\endproclaim

\demo{Proof}  According to \cite{\Hull}, it suffices to
show that all two-sided artinian, left pure semisimple
rings are as claimed.  So suppose that
$R$ is left pure semisimple, artinian on both sides, and
satisfies our product condition.  Since every left
$R$-module and thus, in particular, every direct product of
modules is a direct sum of finitely generated components in
this situation, we infer that any family of pairwise
non-isomorphic finitely generated indecomposable left
$R$-modules is finite.
\qed
\enddemo

	 Next we turn to Artin algebras.  In this setting, our
product property is an immediate consequence of a result
due to Auslander \cite{\Auslarge}.

\proclaim{Proposition 2}  If  $R$ is an Artin algebra and
$(M_n)_{n \in \NN}$ a family of pairwise non-isomorphic
finitely generated, indecomposable $R$-modules, then
$\prod_{n \in \NN} M_n$ is not a direct sum of finitely
generated submodules.
\endproclaim

\demo{Proof} First we observe the following elementary
fact:  If
$R$ is an Artin algebra, and $M = \prod_{n \in \NN} M_n$ a
countable direct product of nonzero left $R$-modules, then
$M$ is not countably generated:  Clearly, this is true when
$R$ is a field, and hence also when $R$ is commutative. 
Suppose that
$R$ is non-commutative, let $C$ be the center of $R$, and
$J(C)$ the Jacobson radical of $C$.  Since $C/J(C)$ is a
finite direct product of fields and $J(C) M_n$ is properly
contained in $M_n$ for each $n$, the $R/J(C)$-module $M /
J(C)M  =  \prod_{n \in \NN} M_n / J(C)M_n$ is not countably
generated by the opening remark.  Hence, neither is $M$ over
$R$.
	Now suppose $(M_n)_{n \in \NN}$ is a family of pairwise
non-isomorphic finitely generated indecomposable modules
over an Artin algebra $R$.  If we had $\prod_{n \in \NN}
M_n =
\bigoplus_{i
\in I} Q_i$ ($\dagger$) with finitely generated $Q_i$, the
above fact would force $I$ to be uncountable.  On the other
hand, each of the $Q_i$ would be isomorphic to some $M_k$ by
Auslander's Cor\. 3.2 \cite{\Auslarge}, whence at least one
$Q_j$ would occur twice in the right-hand side of
($\dagger$), i.e. $Q_j \cong Q_k \cong M_{n(j)}$ with $j
\ne k$.  Now modules with local endomorphism rings can be
canceled from direct sum decompositions, which shows that
$\prod_{n \in \NN \setminus\{n(j)\}} M_n
\cong \bigoplus_{i \in I\setminus\{j\}} Q_i$.  But the
appearance of $Q_j$ in the latter sum contradicts
Auslander's result, in view of the fact that the $M_n$ are
pairwise non-isomorphic. \qed \enddemo

Finally, we prove the product condition for a substantial
class of commutative noetherian rings.

\proclaim{Proposition 3}  Suppose that $R$ is a commutative
noetherian domain of Krull dimension 1 and $(M_n)_{n \in
\NN}$ a family of pairwise non-isomorphic finitely
generated, indecomposable $R$-modules.  Then $\prod_{n \in
\NN} M_n$ is not a direct sum of finitely generated
submodules. \endproclaim

\demo{Proof} Let $R$ and $(M_n)$ be as in the hypothesis,
and assume that $M :=
\prod_{n \in \NN} M_n = \bigoplus_{i \in I} Q_i$ with $Q_i$
finitely generated.  Since, for any nonzero ideal $A$, the
factor ring $R/A$ is artinian, the setup clearly implies
that
$A$ does not annihilate a cofinite subfamily of the family
$(M_n)$; for if $A M_n = 0$ for all but finitely many $n \in
\NN$, say for $n \in N'$, one is dealing with a
decomposition of the module $M/AM = \prod_{n \in N'} M_n
\oplus {\prod_{\NN
\setminus N'}} M_n/AM_n$ over $R/A$ and, after canceling,
one obtains a contradiction to the fact that commutative
artinian rings satisfy the conjecture.  

All $M_n$ being
noetherian, we can find finitely generated torsionfree submodules
$U_n
\subseteq M_n$ such that each of the quotients $M_n / U_n$ is
a torsion module.  Moreover, we let 
$A_n$ be the annihilator of 
$M_n / U_n$ for $n \in \NN$, pick any maximal ideal
$B$ of $R$, and let
$C_1, C_2, C_3, \dots $ be an enumeration of all finite
products of ideals in the set
$\{A_n \mid n \in \NN\}\cup \{B^n \mid n \in \NN\}$. Then
$\bigcap_{n \in \NN} (C_nM_k) = 0$ for all $k \in \NN$, and
thus
$\bigcap_{n \in \NN} (C_nM) = 0$, while, for each $m \in
\NN$,
$\bigcap_{n \le m} (C_n M_k) \ne 0$ for infinitely many $k
\in
\NN$ by the initial comment.  On the other hand, an
extension of Chase's result in  \cite{\ChaII} (see
\cite{\Huibiel, Lemma 11} for a precise statement and proof
of the upgraded version) yields a natural number $L$ such
that
$$\prod_{k \ge L} \bigl( \bigcap_{n \le L} (C_n M_k) \bigr)
\subseteq
\bigoplus_{i \in I'} Q_i + \bigoplus_{i \in I} \bigl(
\bigcap_{n
\in \NN}(C_n Q_i) \bigr),$$
where $I'$ is a finite subset of
$I$.  But the last summand on the right-hand side of this
inclusion is zero by construction, whereas the left-hand
side is still an infinite product.  This forces an infinite
product of nonzero modules over the noetherian ring $R$ to
be finitely generated, an obvious impossibility. \qed 
\enddemo

\noindent {\it Concluding Remarks:} A slight modification
of this argument actually goes through for rings which are
module-finite over a central noetherian subdomain of Krull
dimension 1.  More general results (not just for
commutative rings) should be available, via noetherian
induction, for classes of noetherian rings $R$ having the
following two properties: (i) each proper factor ring of
$R$ satisfies the conjecture, and (ii) for every finitely
generated left $R$ module $M$ there exists a countable
family $(C_n)_{n \in \NN}$ of nonzero ideals such that
$\bigcap_{n \in \NN} (C_n M) = 0$.

\enddocument